\documentclass[twoside]{article} 
\date{} 
\title{Asymptotics of a Mathieu-Gaussian series}
\author{\sc R. B.\ Paris \\
{\em Division of Computing and Mathematics,} \\
{\em Abertay University, Dundee DD1 1HG, UK}}
\begin{document}
\def\f#1#2{\mbox{${\textstyle \frac{#1}{#2}}$}}
\def\dfrac#1#2{\displaystyle{\frac{#1}{#2}}}
\def\boldal{\mbox{\boldmath $\alpha$}}
\newcommand{\bee}{\begin{equation}}
\newcommand{\ee}{\end{equation}}
\newcommand{\lam}{\lambda}
\newcommand{\ka}{\kappa}
\newcommand{\al}{\alpha}
\newcommand{\la}{\lambda}
\newcommand{\ga}{\gamma}
\newcommand{\eps}{\epsilon}
\newcommand{\fr}{\frac{1}{2}}
\newcommand{\fs}{\f{1}{2}}
\newcommand{\g}{\Gamma}
\newcommand{\br}{\biggr}
\newcommand{\bl}{\biggl}
\newcommand{\ra}{\rightarrow}
\newcommand{\gtwid}{\raisebox{-.8ex}{\mbox{$\stackrel{\textstyle >}{\sim}$}}}
\newcommand{\ltwid}{\raisebox{-.8ex}{\mbox{$\stackrel{\textstyle <}{\sim}$}}}
\renewcommand{\topfraction}{0.9}
\renewcommand{\bottomfraction}{0.9}
\renewcommand{\textfraction}{0.05}
\newcommand{\mcol}{\multicolumn}
\date{}
\maketitle
\pagestyle{myheadings}
\markboth{\hfill \sc R. B.\ Paris  \hfill}
{\hfill \sc Asymptotics of a Mathieu-Gaussian series\hfill}
\begin{abstract}
We consider the asymptotic expansion of the functional series 
\[S_{\mu,\ga}(a;\la)=\sum_{n=1}^\infty \frac{n^\ga e^{-\la n^2/a^2}}{(n^2+a^2)^\mu}\]
for real values of the parameters $\ga$, $\la>0$ and $\mu\geq0$ as $|a|\to \infty$ in the sector $|\arg\,a|<\pi/4$. For general values of $\ga$ the expansion is of algebraic type with terms involving the Riemann zeta function and a terminating confluent hypergeometric function. Of principal interest in this study is the case corresponding to even integer values of $\ga$, where the algebraic-type expansion consists of a finite number of terms together with a contribution comprising an infinite sequence of increasingly subdominant exponentially small expansions. This situation is analogous to the well-known Poisson-Jacobi formula corresponding to the case $\mu=\ga=0$. Numerical examples are provided to illustrate the accuracy of these expansions.

\vspace{0.3cm}

\noindent {\bf Mathematics subject classification (2010):} 30E15, 33E20, 33C15, 34E05, 41A60
\vspace{0.1cm}
 
\noindent {\bf Keywords:} Mathieu series, asymptotic expansions, exponential asymptotics, Mellin transform method
\end{abstract}

\vspace{0.3cm}

\noindent $\,$\hrulefill $\,$

\vspace{0.3cm}
\begin{center}
{\bf 1.\ Introduction}
\end{center}
\setcounter{section}{1}
\setcounter{equation}{0}
\renewcommand{\theequation}{\arabic{section}.\arabic{equation}}
The functional series
\bee\label{e11}
\sum_{n=1}^\infty\frac{n}{(n^2+a^2)^\mu}
\ee
in the case $\mu=2$ 
was introduced by Mathieu in his 1890 book \cite{M} dealing with the elasticity of solid bodies. 
The asymptotic expansion for large $a$ of more general functional series has been discussed in \cite{P} and \cite{Z}. 
More recently, Gerhold and Tomovski \cite{G} extended the asymptotic study of such Mathieu series by introducing in (\ref{e11}) the factor $z^n$, where $|z|\leq 1$. From this result they were able to deduce, in particular, the large-$a$ expansions of the trigonometric Mathieu series
\[\sum_{n=1}^\infty \frac{n \sin nx}{(n^2+a^2)^\mu},\qquad\sum_{n=1}^\infty \frac{n \cos nx}{(n^2+a^2)^\mu}.\]
Subsequently, the above trigonometric series were generalised to include the oscillatory Bessel functions $J_\nu(x)$ and $Y_\nu(x)$ with argument proportional to $n/a$, and their large-$a$ asymptotics determined in \cite{P3}. In addition, this last study also considered the inclusion of the modified Bessel function $K_\nu(x)$ of similar argument, which contains the decaying exponential as a special case.  

The asymptotic expansion we consider in this paper is the Mathieu series coupled with a Gaussian exponential of the form
\bee\label{e12}
S_{\mu,\gamma}(a;\lambda):=\sum_{n=1}^\infty\frac{n^\gamma e^{-\la n^2/a^2}}{(n^2+a^2)^\mu}\qquad (\mu\geq0,\ \lambda>0)
\ee
for $|a|\ra\infty$ in the sector $|\arg\,a|<\pi/4$. It will be supposed throughout that $\ga$ is real, although the analysis is easily modified to incorporate complex $\ga$.
We shall employ the Mellin transform approach used in \cite{P,P3,Z}, 
where our interest will be primarily concerned with even integer values of $\gamma$ (positive or negative). We shall find that the asymptotic expansion of $S_{\mu,\gamma}(a;\lambda)$ with these parameter values for large complex $a$ in the sector $|\arg\,a|<\pi/4$  consists of a finite algebraic expansion together with an infinite sequence of increasingly subdominant exponentially small contributions. 

It is interesting that the apparently simple series (\ref{e12}) should possess such an intricate
asymptotic structure in the case of even integer values of $\ga$. This is also found to be the case when $\la=0$ in (\ref{e12}); see \cite{P2} for details. A well-known related series corresponding to $\mu=\ga=0$ is the Poisson-Jacobi formula \cite[p.~124]{WW}
\bee\label{e12a}
S_{0,0}(a;\la)=\sum_{n=1}^\infty e^{-\la n^2/a^2}=\frac{a}{2}\sqrt{\frac{\pi}{\la}}-\frac{1}{2}+a\sqrt{\frac{\pi}{\la}} \sum_{n=1}^\infty e^{-\pi^2n^2a^2/\la}.
\ee
This sum is also seen to consist of a finite algebraic contribution together with an infinite sum of exponentially small terms when $|a|\to\infty$ in the sector $|\arg\,a|<\pi/4$.

In the application of the Mellin transform method to the series in (\ref{e12}) and its alternating variant we shall require the following estimates for the gamma function and the Riemann zeta function. For real $\sigma$ and $t$, we have the estimates
\bee\label{e13}
\g(\sigma\pm it)=O(t^{\sigma-\frac{1}{2}}e^{-\frac{1}{2}\pi t}),\qquad |\zeta(\sigma\pm it)|=O(t^{\Omega(\sigma)} \log^\alpha t)\quad (t\ra+\infty),
\ee
where $\Omega(\sigma)=0$ ($\sigma>1$), $\fs-\fs\sigma$ ($0\leq\sigma\leq 1$), $\fs-\sigma$ ($\sigma<0$) and $\alpha=1$ ($0\leq\sigma\leq 1$), $\alpha=0$ otherwise \cite[p.~95]{ECT}. The zeta function $\zeta(s)$ has a simple pole of unit residue at $s=1$ and the evaluations for positive integer $k$
\bee\label{e14}
\zeta(0)=-\fs,\quad \zeta(-2k)=0,\quad \zeta(2k)=\frac{(2\pi)^{2k}}{2(2k)!}|B_{2k}|\ \ (k\geq 1),
\ee
where $B_k$ are the Bernoulli numbers. Finally, we have the well-known functional relation satisfied by $\zeta(s)$ given by \cite[p.~603]{DLMF}
\bee\label{e15}
\zeta(s)=2^s \pi^{s-1} \zeta(1-s) \g(1-s) \sin \fs\pi s.
\ee
\vspace{0.6cm}

\begin{center}
{\bf 2.\ An integral representation}
\end{center}
\setcounter{section}{2}
\setcounter{equation}{0}
\renewcommand{\theequation}{\arabic{section}.\arabic{equation}}
The generalised Mathieu series defined in (\ref{e12}) can be written as
\bee\label{e21}
S_{\mu,\gamma}(a;\lambda)=a^{-\delta} \sum_{n=1}^\infty h(n/a),\qquad h(x):=\frac{x^\gamma e^{-\lambda x^2}}{(1+x^2)^\mu},\quad \delta:=2\mu-\gamma,
\ee
where the parameter $\delta$ is real and $\la>0$. We employ a Mellin transform approach as discussed in \cite[Section 4.1.1]{PK}. The Mellin transform of $h(x)$ is ${\cal H}(s)=\int_0^\infty x^{s-1}h(x)\,dx$, where \cite[(13.4.4)]{DLMF}
\bee\label{e21a}
{\cal H}(s)=\int_0^\infty \frac{x^{\gamma+s-1}e^{-\la x^2}}{(1+x^2)^\mu}\,dx=\fs\g(\f{\ga+s}{2})\,U(\f{\ga+s}{2},1+\f{\ga+s}{2}-\mu,\la)
\ee
in the half-plane $\Re (s)>-\ga$,
with $U(a,b,z)$ being the confluent hypergeometric function of the second kind. The transform ${\cal H}(s)$ can be represented alternatively in the form
\bee\label{e21b}
{\cal H}(s)=\fs\{{\cal H}_1(s)+{\cal H}_2(s)\},
\ee
where
\bee\label{e21c}
{\cal H}_1(s)=\frac{\g(\frac{\ga+s}{2}) \g(\mu-\frac{\ga+s}{2})}{\g(\mu)}\,{}_1F_1(\f{\ga+s}{2};1+\f{\ga+s}{2}-\mu;\la),
\ee
\bee\label{e21d}
{\cal H}_2(s)=\la^{\mu-(\ga+s)/2}\g(\f{\ga+s}{2}-\mu)\,{}_1F_1(\mu;1-\f{\ga+s}{2}+\mu;\la).
\ee

Using the Mellin inversion theorem (see, for example, \cite[p.~118]{PK}), we find
\bee\label{e22}
S_{\mu,\ga}(a;\la)=\frac{a^{-\delta}}{2\pi i}\sum_{n=1}^\infty \int_{c-\infty i}^{c+\infty i} {\cal H}(s) (n/a)^{-s}ds=\frac{a^{-\delta}}{2\pi i}\int_{c-\infty i}^{c+\infty i} {\cal H}(s) \zeta(s) a^s ds,
\ee
where $\zeta(s)$ is the Riemann zeta function and $c>\max\{1,-\ga\}$. The inversion of the order of summation and integration is justified by absolute convergence provided $c$ satisfies this condition. 
From the estimates in (\ref{e13}) and the fact that the confluent hypergeometric functions\footnote{The ${}_1F_1$ function appearing in ${\cal H}_1(s)$ can be written as $e^\la\,{}_1F_1(1-\mu;1+\fs(\ga+s)-\mu;-\la)$ by Kummer's transformation \cite[p.~325]{DLMF}.}  appearing in ${\cal H}_1(s)$ and ${\cal H}_2(s)$ are both $O(1)$ as $\Im(s)\to\pm\infty$, the integral in (\ref{e22}) then defines $S_{\mu,\gamma}(a;\lambda)$ for complex $a$ in the sector $|\arg\,a|<\pi/4$. 
The integration path in (\ref{e22}) lies to the right of the simple pole of $\zeta(s)$ at $s=1$ and the poles of $\g((\ga+s)/2)$ at $s=-\gamma-2k$ ($k=0, 1, 2, \ldots\,$), these being the only poles of the integrand since the $U$ function in (\ref{e21a}) has no poles; see Appendix A for a demonstration of this fact.  
 
We consider the integral in (\ref{e22}) taken round the rectangular contour with vertices at $c\pm iT$ and $-c'\pm iT$,
where $c'>0$ and $T>0$. The contribution from the upper and lower sides of the rectangle $s=\sigma\pm iT$, $-c'\leq\sigma\leq c$, vanishes as $T\ra\infty$ provided $|\arg\,a|<\pi/4$, since from (\ref{e13}), the modulus of the integrand is controlled by $O(T^{\Omega(\sigma)+(\sigma-\delta-1)/2} \log\,T\,e^{-\Delta T})$, where $\Delta=\pi/4-|\arg\,a|$. Displacement of the integration path to the left over the pole at $s=1$ and those of ${\cal H}_1(s)$ at $s=-2k-\ga$ (when $\mu>0$) then yields
\[S_{\mu,\gamma}(a;\lambda)- a^{1-\delta} {\cal H}(1)\hspace{8cm}\]
\bee\label{e24}
\hspace{4cm}\sim a^{-2\mu}\sum_{k=0}^\infty \frac{(-)^k (\mu)_k}{a^{2k}} \zeta(-2k-\ga)\,{}_1F_1(-k;1-\mu-k;\la),
\ee
where
\bee\label{e26}
{\cal H}(1)=\frac{1}{2}\g(\f{1+\ga}{2})\,U(\f{1+\ga}{2},1+\f{1+\ga}{2}-\mu,\la).
\ee
When $\mu=0$, we have ${\cal H}_1(s)\equiv 0$ and the poles
from ${\cal H}_2(s)$ at $s=-2k-\ga$ yield
\bee\label{e25}
S_{0,\ga}(a;\la)-\frac{1}{2}\bl(\frac{a^2}{\la}\br)^{(\ga+1)/2}\,\g(\frac{\ga+1}{2})\sim\sum_{k=0}^\infty \frac{(-\la)^k}{k!\,a^{2k}}\,\zeta(-2k-\ga),
\ee
where we have used the fact that $U(\al,1+\al,x)=x^{-\al}$.
Both expansions (\ref{e24}) and (\ref{e25}) hold as $|a|\to\infty$ in $|\arg\,a|<\pi/4$. 
The hypergeometric functions appearing in the sum in (\ref{e24}) are polynomials in $\la$ of degree $k$ when $\mu$ is non-integer; for positive integer $\mu$ they can be expressed by Kummer's theorem as $e^\la$ multiplied by a polynomial in $\la$ of degree $\mu-1$.

The expansion (\ref{e24}) holds for general values of $\mu>0$ and $\gamma$. If $\ga$ equals an odd negative integer, the pole at $s=1$ has both a double pole contribution (resulting from ${\cal H}_1(s)$) and a simple pole contribution (resulting from ${\cal H}_2(s)$). An example of the expansion when $\ga=-1$ is discussed in Appendix B. We remark that when $\ga=0$ in (\ref{e25}) the expansion for $S_{0,0}(a;\la)$ correctly reduces to the first two terms in the Poisson-Jacobi formula (\ref{e12a}), but does not account for the exponentially small contribution. Further consideration of this case is discussed at the end of Section 3.

Finally, we note that when $\ga=2p$ is an even integer, there will be a finite number of poles of the integrand of the sequence $s=-2k-2p$ on account of the trivial zeros of $\zeta(s)$ at $s=-2, -4, \ldots\ $.
This results in the number of terms in the asymptotic series in (\ref{e24}) and (\ref{e25}) being either finite or zero. This situation is the main subject of this paper. 
We shall show that, in addition to a finite algebraic contribution, there is a sequence of increasingly subdominant exponentially small terms in the large-$a$ limit. This is analogous to the exponentially small contribution appearing on the right-hand side of the Poisson-Jacobi formula (\ref{e12a}).

\vspace{0.6cm}

\begin{center}
{\bf 3.\ The exponentially small contribution to $S_{\mu,\ga}(a;\la)$ when $\ga=2p$}
\end{center}
\setcounter{section}{3}
\setcounter{equation}{0}
\renewcommand{\theequation}{\arabic{section}.\arabic{equation}}
Let $\ga=2p$ be an even integer and $\mu\geq0$. Then the quantity $\delta$ defined in (\ref{e21}) is $\delta=2(\mu-p)$. The number of poles of the sequence $s=-2k-2p$ is finite (when $p\leq 0$) or zero (when $p\geq 1$)  on account of the trivial zeros of $\zeta(s)$. Then we have upon displacement of the integration path
\bee\label{e31}
S_{\mu,\ga}(a;\la)=a^{1-\delta} {\cal H}(1)+H_{\mu,\ga}(a;\la)+J(a;\la),
\ee
where
\bee\label{e31a}
H_{\mu,\ga}(a;\la)=
\left\{\begin{array}{ll} 0 & (p\geq1)\\
-\fs a^{-\delta} & (p=0)\\
a^{-\delta}\sum_{k=0}^{|p|} \zeta(2k)R_k(\mu,|p|) a^{2k} & (p\leq -1).\end{array}\right.
\ee
The quantity $R_k(\mu,|p|)$ denotes the residue of ${\cal H}(s)$ at $s=2k$, $0\leq k\leq |p|$ when $p\leq -1$ given by
\[R_k(\mu,|p|)=\frac{(-)^{q-k}}{(q-k)!}\,U(k-q,1+k-q-\mu,\la) \qquad (q=-p,\ 0\leq k\leq q).\]
Routine calculations show that when $p=-1, -2$, for example, we have
\begin{eqnarray*}
R_0(\mu,1)\!\!&=&\!\!-(\mu+\la),\quad R_1(\mu,1)=1,\\
R_0(\mu,2)\!\!&=&\!\!\fs\mu(1+\mu)+\mu\la+\fs\la^2,\quad R_1(\mu,2)=-(\mu+\la), \quad R_2(\mu,2)=1.
\end{eqnarray*}

The integral $J(a;\la)$ appearing in (\ref{e31}) is defined by
\[J(a;\la)=\frac{a^{-\delta}}{2\pi i} \int_{c-\infty i}^{c+\infty i} {\cal H}(-s) \zeta(-s)a^{-s}ds\qquad (c>0),\]
where we have replaced $s$ by $-s$. Use of the functional relation for $\zeta(s)$ in (\ref{e15}) followed by expansion of $\zeta(1+s)$ (permissible since $c>0$) leads to
\begin{eqnarray*}
J(a;\la)&=&-\frac{a^{-\delta}}{2\pi i}\int_{c-\infty i}^{c+\infty i} {\cal H}(-s)\zeta(1+s)\g(1+s) \frac{\sin \fs\pi s}{\pi}\,(2\pi a)^{-s}ds\nonumber\\
&=&-\frac{a^{-\delta}}{\pi} \sum_{k\geq 1} \frac{1}{k}\cdot \frac{1}{2\pi i}\int_{c-\infty i}^{c+\infty i} {\cal H}(-s)\g(1+s)\sin \fs\pi s\, (2\pi ka)^{-s}ds.
\end{eqnarray*}
From (\ref{e21a}) and an application of Kummer's transformation \cite[(13.2.40)]{DLMF} we have 
\[{\cal H}(-s)=\fs \g(p-\fs s)\,\la^{\mu-p+s/2} U(\mu,1+\mu-p+\fs s,\la),\]
whence
\bee\label{e32}
J(a;\la)=\frac{(-)^pa^{-\delta}}{2\la^{p-\mu}} \sum_{k\geq 1} \frac{1}{k}\cdot \frac{1}{2\pi i}\int_{c-\infty i}^{c+\infty i} \frac{2^{-s}\g(1+s)}{\g(1-p+\fs s)}\,X_k^{-s/2}\,U(\mu,1\!+\!\mu\!-\!p\!+\!\fs s,\la)\,ds,
\ee
where
\bee\label{e32X}
X_k:=\frac{\pi^2k^2a^2}{\la}.
\ee

We first consider the case $p=0, 1, 2, \ldots\,$. Since there are no poles of the integrand in (\ref{e32}) in $\Re (s)>0$ the integration path can be displaced as far to the right as we please such that $|s|$ is everywhere large on the new path. The quotient of gamma functions in the integrand can then be expanded as \cite[p.~53]{PK}, \cite[(5.11.19)]{DLMF}
\bee\label{e33d}\frac{2^{-s}\sqrt{\pi} \g(1+s)}{\g(1\!-\!p\!+\!\fs s)}=\frac{\g(\fs+\fs s)\g(1+\fs s)}{\g(1\!-\!p\!+\!\fs s)}=\sum_{j=0}^p(-)^j c_j \g(\fs s+\vartheta-j),
\ee
where $\vartheta=p+\fs$. The coefficients are given explicitly by
\bee\label{e33a}
c_j=\frac{1}{j!} (-p)_j (-p+\fs)_j=\frac{(-2p)_{2j}}{2^{2j} j!}.
\ee
We observe that $c_j=0$ for $j>p$ so that the above sum of gamma functions terminates and so is {\it exact}.
Substitution of the expansion (\ref{e33d}) in (\ref{e32}), combined with the integral representation \cite[(13.4.4)]{DLMF}
\[U(\mu,1+\mu-p+\fs s,\la)=\frac{1}{\g(\mu)}\int_0^\infty e^{-\la t}t^{\mu-1} (1+t)^{-p-s/2}dt \qquad (\mu>0),\]
then shows that, provided $\mu>0$,
\[J(a;\la)=\frac{(-)^pa^{-\delta}\la^{p-\mu}}{2\sqrt{\pi} \g(\mu)} \sum_{k\geq 1}\frac{1}{k}\sum_{j=0}^p(-)^j c_j\hspace{5cm}\] 
\[\hspace{3cm}\times\int_0^\infty \frac{e^{-\la t}t^{\mu-1}}{(1+t)^p} \bl(\frac{1}{2\pi i}\int_{c-\infty i}^{c+\infty i}\g(\fs s+\vartheta-j) \bl(\frac{X_k}{1+t}\br)^{\!-s/2}ds\br)dt.\]

The inner integral appearing in $J(a;\la)$ can be evaluated by making use of the well-known result 
\bee\label{e34aa}
\frac{1}{2\pi i}\int_{L} \g(s+\alpha) z^{-s}ds=z^\alpha e^{-z}\qquad (|\arg\,z|<\fs\pi),
\ee
where $L$ is a path parallel to the imaginary $s$-axis lying to the right of all the poles of $\g(s+\alpha)$;
see, for example, \cite[Section 3.3.1]{PK}. Evaluation of the inner integral when $|\arg\,a|<\pi/4$ then produces the final {\it exact\/} result
\[J(a;\la)=\frac{(-)^p\pi^\mu e^\la}{\g(\mu)} \bl(\frac{\la}{\pi a^2}\br)^{\!\mu\!-\!2p\!-\!1/2}\sum_{k\geq 1}k^{2p}e^{-2\pi ka} \sum_{j=0}^p(-)^jc_j\bl(\frac{\la}{\pi^2k^2a^2}\br)^{\!j}I_{jk},\]
where
\bee\label{e320}
I_{jk}=\int_0^\infty \frac{t^{\mu-1}e^{-\psi(t)}}{(1+t)^{2p-j+1/2}}\,dt,\qquad\psi(t):=\la(1+t)+\frac{X_k}{1+t}-2\pi ka.
\ee

In the special case $\mu=0$, we have from (\ref{e32}) (since $U(0,b,z)=1$) that
\[J(a;\la)=\frac{(-)^pa^{2p}}{2\sqrt{\pi}\la^{p}} \sum_{k\geq 1} \frac{1}{k}\sum_{j=0}^p(-)^j c_j\cdot \frac{1}{2\pi i}\int_{c-\infty i}^{c+\infty i} \g(\fs+\vartheta-j)\,X_k^{-s/2}\,ds.\]
Evaluation of the integral by means of (\ref{e34aa}) then produces
\[
J(a;\la)=(-)^p \bl(\frac{\la}{\pi a^2}\br)^{\!-2p-1/2}\sum_{k\geq 1} k^{2p} e^{-\pi^2k^2a^2/\la}\sum_{j=0}^p(-)^j c_j\bl(\frac{\la}{\pi^2k^2a^2}\br)^{\!j} \qquad (\mu=0).\]
Then we have the following theorem:
\newtheorem{theorem}{Theorem}
\begin{theorem}$\!\!\!.$ For $\mu\geq0$, $\la>0$, $\delta=2\mu-\ga$ and $\ga=2p$, where $p$ is a non-negative integer, we have when $|\arg\,a|<\pi/4$
\[S_{\mu,\ga}(a;\la)=a^{1-\delta} {\cal H}(1)-\fs a^{-\delta}\delta_{0p}+J(a;\la),\]
where ${\cal H}(1)$ is defined in (\ref{e26}) and $\delta_{0p}$ is the Kronecker symbol. The exponentially small contribution $J(a;\la)$ is given exactly by the double sums
\bee\label{e321}
J(a;\la)=\frac{(-)^p\pi^\mu e^\la}{\g(\mu)} \bl(\frac{\la}{\pi a^2}\br)^{\!\mu\!-\!2p\!-\!1/2}\sum_{k\geq 1}k^{2p}e^{-2\pi ka} \sum_{j=0}^p(-)^jc_j\bl(\frac{\la}{\pi^2k^2a^2}\br)^{\!j}I_{jk}\qquad (\mu>0),
\ee
and
\bee\label{e322}
J(a;\la)=(-)^p \bl(\frac{\la}{\pi a^2}\br)^{\!-2p-1/2}\sum_{k\geq 1} k^{2p} e^{-\pi^2k^2a^2/\la}\sum_{j=0}^p(-)^j c_j\bl(\frac{\la}{\pi^2k^2a^2}\br)^{\!j} \qquad (\mu=0),
\ee
where the coefficients $c_j=(-2p)_j/(2^{2j}j!)$ and the integrals $I_{jk}$ are defined in (\ref{e320}).
\end{theorem}
\medskip

\noindent{\bf Remark 1.}\ \ When $\mu=p=0$, we find from (\ref{e26}), (\ref{e31}) and (\ref{e322}) (since $U(\fs,\f{3}{2},\la)=\la^{-1/2}$) the result
\[S_{0,0}(a;\la)= \frac{a}{2}\sqrt{\frac{\pi}{\la}}-\frac{1}{2}+a\sqrt{\frac{\pi}{\la}} \sum_{k\geq 1} e^{-\pi^2k^2a^2/\la},\]
which is the Poisson-Jacobi formula stated in (\ref{e12a}). 
\medskip

When $p=1$, we have $U(\f{3}{2},\f{5}{2},\la)=\la^{-3/2}$ and
\[S_{0,2}(a;\la)=\frac{a^3\sqrt{\pi}}{4\la^{3/2}}-\bl(\frac{\pi a^2}{\la}\br)^{\!5/2} \sum_{k\geq 1} \bl(k^2-\frac{\la}{2\pi^2a^2}\br) e^{-\pi^2k^2a^2/\la}.\]
We observe that this last case can also be obtained by differentiation of the Poisson-Jacobi formula with respect to $\la$, since
\[S_{0,2p}(a;\la)=(-)^p a^{2p}\,\frac{\partial^p}{\partial\la^p} S_{0,0}(a;\la) \qquad (p\geq 1).\]
\medskip

\noindent{\bf Remark 2.}\ \ When $p=-1, -2, \ldots\ $, the expansion (\ref{e33d}) does not terminate and becomes an inverse factorial expansion. Then we have the exponentially small contribution given by
\bee\label{e323}
J(a;\la)\sim\frac{(-)^p\pi^\mu e^\la}{\g(\mu)} \bl(\frac{\la}{\pi a^2}\br)^{\!\mu\!-\!2p\!-\!1/2}\sum_{k\geq 1}k^{2p}e^{-2\pi ka} \sum_{j=0}^\infty(-)^jc_j\bl(\frac{\la}{\pi^2k^2a^2}\br)^{\!j}I_{jk}\qquad (\mu>0)
\ee
as $|a|\to\infty$ in $|\arg\,a|<\pi/4$.
\vspace{0.6cm}

\begin{center}
{\bf 4.\ Alternative form of expansion for positive integer values of $\mu$}
\end{center}
\setcounter{section}{4}
\setcounter{equation}{0}
\renewcommand{\theequation}{\arabic{section}.\arabic{equation}}
Let $\mu=m$ be a positive integer and $\ga=2p$.
We split ${\cal H}(s)$ into its two constituent parts given by (\ref{e21b}) and write $J(a;\la)=J_1(a;\la)+J_2(a;\la)$ in an obvious manner. 
\medskip

\noindent {\bf {4.1}\ Evaluation of $J_1(a;\la).$}\ \ Then we have
\[J_1(a;\la)=-\frac{a^{-\delta}}{2\pi}\sum_{k\geq 1}\frac{1}{k}\cdot\frac{1}{2\pi i}\int_{c-\infty i}^{c+\infty i} {\cal H}_1(-s)\g(1+s) \sin \fs\pi s\,(2\pi ka)^{-s}ds,\] 
where from (\ref{e21c}) with an application of Kummer's transformation
\begin{eqnarray*}
{\cal H}_1(-s)&=&\frac{(-)^{m-p}\pi e^\la}{\g(m)\,\sin \fs\pi s}\,\frac{\g(p-\fs s)}{\g(1-m+p-\fs s)}\,{}_1 F_1(1-m;1-m+p-\fs s;-\la)\\
&=&\frac{(-)^{m-p}\pi e^\la}{\g(m)\,\sin \fs\pi s}\,\sum_{r=0}^{m-1}\frac{(1-m)_r (-\la)^r}{r!}\,\frac{\g(p-\fs s)}{\g(1-m+p-\fs s+r)}\\
&=&\frac{(-)^{p-1}\pi e^\la}{\g(m)\,\sin \fs\pi s}\,\sum_{r=0}^{m-1}\frac{(1-m)_r \la^r}{r!}\,\frac{\g(m-p-r+\fs s)}{\g(1-p+\fs s)}.
\end{eqnarray*}
If we make the change of summation index $r\to m-1-\ell$ and use the fact that $(1-m)_{m-1-\ell}=(-)^{m-1+\ell} \g(m)/\ell!$, we find
\bee\label{e33}
{\cal H}_1(-s)=\frac{(-)^{m-p}\pi e^\la}{\sin \fs\pi s}\,\sum_{\ell=0}^{m-1} A_\ell (1-p+\fs s)_\ell,\qquad A_\ell:=\frac{(-)^\ell \la^{m-1-\ell}}{(m-1-\ell)! \ell!}.
\ee
The Pochhammer symbol appearing in (\ref{e33}) can be written in the form
\[(1-p+\fs s)_\ell=\sum_{r=0}^\ell B_{r\ell} (s+1)_r,\]
where
\[B_{00}=1,\quad B_{01}=\fs-p,\quad B_{11}=\fs\]
\[B_{02}=(\fs-p)_2,\quad B_{12}=\f{3}{4}-p,\quad B_{22}=\f{1}{4},\]
\bee\label{e33c}
B_{03}=(\fs-p)_3,\quad B_{13}=\f{3}{8}(5-10p+4p^2),\quad B_{23}=\f{3}{4}(1-p), \quad B_{33}=\f{1}{8},\  \ldots\ .
\ee

Then we obtain
\bee\label{e33ab}
J_1(a;\la)=(-)^{k-p-1}\frac{e^\la a^{-\delta}}{2} \sum_{\ell=0}^{m-1}\sum_{r=0}^\ell A_\ell B_{r\ell} \sum_{k\geq 1}\frac{1}{k}\cdot \frac{1}{2\pi i}\int_{c-\infty i}^{c+\infty i}\g(1+s+r)(2\pi ka)^{-s}ds.
\ee 
The integrals appearing in $J_1(a;\la)$ can be evaluated by (\ref{e34aa})
to produce the final {\it exact\/} result
\bee\label{e34}
J_1(a;\la)=(-)^{m-p-1}a^{2p+1}\pi e^\la  \sum_{r=0}^{m-1}\sum_{\ell=r}^{m-1} A_\ell B_{r\ell}\frac{\sigma_r(a)(2\pi
 a)^r}{a^{2m-r}}
\ee
for positive integer $m$. Here we have defined the sums 
\bee\label{e34c}
\sigma_r(a):=e^{2\pi a}\sum_{k\geq 1} k^r e^{-2\pi ka},
\ee
which have the evaluations for $0\leq r\leq 3$
\[\sigma_0(a)=\frac{e^{\pi a}}{2\sinh \pi a},\quad \sigma_1(a)=\frac{e^{2\pi a}}{4 \sinh^2 \pi a}, \quad \sigma_2(a)=\frac{e^{2\pi a}\cosh \pi a}{4\sinh^3 \pi a},\]
\[\sigma_3(a)=\frac{e^{2\pi a}(2+\cosh \pi a)}{8 \sinh^4 \pi a}, \ldots\ .\]
Note that $J_1(a;\la)\equiv 0$ when $m=0$, since ${\cal H}_1(s)$ vanishes for these values.
\bigskip

\noindent {\bf {4.2}\ Evaluation of $J_2(a;\la).$}\ \ We have 
\[J_2(a;\la)=-\frac{a^{-\delta}}{2\pi}\sum_{k\geq 1}\frac{1}{k}\cdot\frac{1}{2\pi i}\int_{c-\infty i}^{c+\infty i} {\cal H}_2(-s)\g(1+s) \sin \fs\pi s\,(2\pi ka)^{-s}ds,\] 
where from (\ref{e21d})
\begin{eqnarray*}
{\cal H}_2(s)&=&(-)^{m-p-1}\frac{\pi \la^{m-p-s/2}}{\sin \fs\pi s}\,\frac{{}_1F_1(m;1+m-p+\fs s;\la)}{\g(1+m-p+\fs s)}\\
&=&(-)^{m-p-1}\frac{\pi \la^{m-p-s/2}}{\sin \fs\pi s}\,\sum_{r=0}^\infty\frac{(m)_r \la^r}{r! \g(1\!+\!m\!-\!p\!+\!r\!+\!\fs s)}.
\end{eqnarray*}
Then we obtain
\bee\label{e35b}
J_2(a;\la)=\frac{(-\la)^{m-p}a^{-\delta}}{\sqrt{\pi}} \sum_{k\geq 1}\frac{1}{k}\sum_{r=0}^\infty\frac{(m)_r \la^r}{r!}\,K_{kr},
\ee
where
\bee\label{e35}
K_{kr}=\frac{\sqrt{\pi}}{4\pi i}\int_{c-\infty i}^{c+\infty i}\frac{2^{-s} \g(1+s)}{\g(1\!+\!m\!-\!p\!+\!r\!+\!\fs s)}\,X_k^{-s/2}ds
\ee
with $X_k$ defined in (\ref{e32X}). 
Since there no poles of the integrand in $\Re (s)>0$ the integration path in (\ref{e35}) can be displaced as far to the right as we please such that $|s|$ is everywhere large on the new path. The quotient of gamma functions in the integrand can then be expanded in a manner similar to that in (\ref{e33d}) to find  \cite[p.~53]{PK}
\[\frac{2^{-s}\sqrt{\pi} \g(1+s)}{\g(1\!+\!m\!-\!p\!+\!r\!+\!\fs s)}=\frac{\g(\fs+\fs s)\g(1+\fs s)}{\g(1\!+\!m\!-\!p\!+\!r\!+\!\fs s)}\hspace{4cm}\]
\[\hspace{3cm}=\sum_{j=0}^{M-1}(-)^j {\hat c}_j(r) \g(\fs+\vartheta-j)+\rho_M(s) \g(\fs s+\vartheta-M),\]
where $M$ is a positive integer, $\vartheta=\fs+p-m-r$ and $\rho_M(s)=O(1)$ as $|s|\to\infty$ in $|\arg\,s|<\pi$. The coefficients are given explicitly by
\bee\label{e36}
{\hat c}_j(r)=\frac{1}{j!} (m-p+r)_j (m-p+r+\fs)_j=\frac{(2m-2p+2r)_{2j}}{2^{2j} j!}.
\ee
When $m=r=0$, the coefficients ${\hat c}_j(r)$ reduce to $c_j$ in (\ref{e33a}).
Then
\[K_{kr}=\sum_{j=0}^{M-1}(-)^r c_j(r)\cdot \frac{1}{4\pi i}\int_{c-\infty i}^{c+\infty i}\g(\fs s+\vartheta-j) X_k^{-s/2}ds+R_{M,r}\]
\bee\label{e36a}
=X_k^\vartheta e^{-X_k}\bl\{\sum_{j=0}^{M-1}(-)^j c_j(r) X_k^{-j}+O(X_k^{-M})\br\}
\ee
by (\ref{e34aa}), where the remainder term
\[R_{M,r}=\frac{1}{4\pi i}\int_{c-\infty i}^{c+\infty i}\rho_M(s) \g(\fs s+\vartheta-j) X_k^{-s/2}ds=O(X_k^{\vartheta-M}e^{-X_k})\]
as $|a|\to\infty$ in $|\arg\,a|<\pi/4$ by Lemma 2.7 in \cite[p.~71]{PK}.

From (\ref{e35b}) and (\ref{e36a}), it then follows that
\[J_2(a;\la)\sim (-)^{m-p}\bl(\frac{\la}{\pi a^2}\br)^{\!\delta-1/2}\sum_{k\geq 1}\frac{e^{-\pi^2k^2a^2/\la}}{k^\delta}\sum_{r=0}^\infty \sum_{j=0}^\infty (-)^j c_j(r)\,\frac{(m)_r \la^r}{r!}\bl(\frac{\la}{\pi^2k^2a^2}\br)^{\!r+j}\]
\bee\label{e3j}
=(-)^{m-p}\bl(\frac{\la}{\pi a^2}\br)^{\!\delta-1/2}\sum_{r=0}^\infty (-)^rC_r\bl(\frac{\la}{\pi^2a^2}\br)^{\!r} \sum_{k\geq 1}\frac{e^{-\pi^2k^2a^2/\la}}{k^{2r+\delta}},
\ee
where the double sum over $r$ and $j$ has been summed `diagonally' (see \cite[p.~58]{S}) and the coefficients $C_r$ are given by
\[C_r\equiv C_r(m,p,\la)=\sum_{n=0}^r\frac{(-\la)^n (m)_n}{n!}\,{\hat c}_{r-n}(n).\]
Use of (\ref{e36}) shows that the $C_r$ may be expressed in terms of a terminating ${}_2F_2$ hypergeometric series:
\[C_r=\sum_{n=0}^r\frac{(-\la)^r (m)_n}{n! (r-n)!}\,\frac{2^{2n-2r}(2m-2p)_{2r}}{(2m-2p)_{2n}}=\frac{(2m-2p)_{2r}}{2^{2r} r!}\sum_{n=0}^r\frac{(-r)_n(m)_n \la^n}{n! (m-p)_n (m-p+\fs)_n}\]
\bee\label{e37}
=\frac{(2m-2p)_{2r}}{2^{2r} r!}\,{}_2F_2(-r,m;m-p,m-p+\fs;\la),
\ee
where we have made use of the results $(a+2n)_{2r-2n}=(a)_{2r}/(a)_{2n}$ and $(2a)_{2n}=2^{2n}(a)_n (a+\fs)_n$.

Then we have the following theorem:
\begin{theorem}$\!\!\!.$
Let $\mu=m$ and $\gamma=2p$, where $m\geq 1$, $p$ are integers, and $\delta=2(m-p)$ with $\la>0$. Then we have the representation
\bee\label{e38a}
S_{\mu,\gamma}(a;\lambda)=a^{1-\delta}{\cal H}(1)+H_{\mu,\ga}(a;\la)+J(a;\la),
\ee
where ${\cal H}(1)$ and $H_{\mu,\ga}(a;\la)$ are defined in (\ref{e26}) and (\ref{e31a}). The exponentially small contribution $J(a;\la)$ has the expansion
\[J(a;\la)+(-)^{m-p}\pi a^{2p+1}e^{\la-2\pi a}\sum_{r=0}^{m-1}\sum_{\ell=r}^{m-1}A_\ell B_{r\ell}\,\frac{\sigma_r(a) (2\pi)^r}{a^{2m-r}}\hspace{4cm}\]
\bee\label{e38b}
\hspace{3cm}\sim(-)^{m-p}\bl(\frac{\la}{\pi a^2}\br)^{\!\delta-1/2}\sum_{r=0}^\infty(-)^r C_r\bl(\frac{\la}{\pi^2 a^2}\br)^{\!r}\sum_{k\geq1}\frac{e^{-\pi^2k^2a^2/\la}}{k^{2r+\delta}}
\ee
as $|a|\to\infty$ in $|\arg\,a|<\pi/4$. The coefficients $A_\ell$, $B_{r\ell}$ and $C_r$ are defined in (\ref{e33}), (\ref{e33c}) and (\ref{e37}) and the functions $\sigma_r(a)$ are given by (\ref{e34c}). When $m=0$ the double sum on the left-hand side of (\ref{e38b}) vanishes.
\end{theorem}

\vspace{0.6cm}

\begin{center}
{\bf 5.\ Examples}
\end{center}
\setcounter{section}{5}
\setcounter{equation}{0}
\renewcommand{\theequation}{\arabic{section}.\arabic{equation}}
We present some examples of the expansion of $S_{\mu,\gamma}(a;\lambda)$ stated in Theorem 2.
\vspace{0.2cm}

\noindent {\bf Example 1.}\ \ Let $\mu=0$ and $\ga=2p$. Then from (\ref{e38a}) and (\ref{e38b}) we have, for $p=1, 2, \ldots\ $,
\[S_{0,2p}(a;\la)=\sum_{n=1}^\infty n^{2p} e^{-\la n^2/a^2}\sim \frac{1}{2}\bl(\frac{\la}{a^2}\br)^{-p-1/2}\g(p+\fs)\hspace{3cm}\]
\bee\label{e550}
+(-)^{p}\bl(\frac{\la}{\pi a^2}\br)^{\!-2p-1/2}\sum_{r=0}^\infty(-)^r C_r\bl(\frac{\la}{\pi^2 a^2}\br)^{\!r}\sum_{k\geq1}\frac{e^{-\pi^2k^2a^2/\la}}{k^{2r-2p}}
\ee
as $|a|\to\infty$ in $|\arg\,a|<\pi/4$,
where from (\ref{e37}) $C_r=(-2p)_{2r}/(2^{2r} r!)$. The case $p=0$ is covered in Remark 1.

In the case $p\leq -1$ we let $p=-q$ and note that ${\cal H}(s)=\fs\la^{q-s/2} \g(\fs s-q)$ when $\mu=0$. The residue of ${\cal H}(s)$ at $s=2q-2k$ ($k=0, 1, 2, \ldots$) is given by $(-)^k/k!$. Then
\[S_{0,-2q}(a;\la)=\sum_{n=1}^\infty \frac{e^{-\la n^2/a^2}}{n^{2q}}\sim \frac{1}{2}\bl(\frac{\la}{ a^2}\br)^{\!q-1/2}\g(\fs-q)+\sum_{k=0}^q\frac{(-)^k}{k!}\bl(\frac{\la}{a^2}\br)^{\!k} \zeta(2q-2k)\]
\bee\label{e43}
+(-)^q\bl(\frac{\la}{\pi a^2}\br)^{\!2q-1/2}\sum_{r=0}^\infty \frac{(-)^r(2q)_{2r}}{2^{2r} r!}\bl(\frac{\la}{\pi^2a^2}\br)^{\!r} \sum_{k\geq 1}\frac{e^{-\pi^2k^2a^2/\la}}{k^{2q+2r}}
\ee
as $|a|\to\infty$ in $|\arg\,a|<\pi/4$.
\bigskip

\noindent {\bf Example 2.}\ \ Let $\mu=1$ and $\ga=0$ (so that $\delta=2$). Then we have
\bee\label{e44}
S_{1,0}(a;\la)=\sum_{n=1}^\infty \frac{e^{-\la n^2/a^2}}{n^2+a^2}=\frac{\pi e^\la}{2a} \mbox{erfc} \sqrt{\la}-\frac{1}{2a^2}+J(a;\la),
\ee
since \cite[(13.6.8)]{DLMF}
\[{\cal H}(1)=\fs\sqrt{\pi}\,U(\fs,\fs,\la)=\fs\pi e^\la\,\mbox{erfc} \sqrt{\la},\]
where erfc is the complementary error function. 
From (\ref{e38b}), the exponentially small contribution is
\bee\label{e45}
J(a;\la)-\frac{\pi e^\la}{2a}\,\frac{e^{-\pi a}}{\sinh \pi a}\sim\bl(\frac{\la}{\pi a^2}\br)^{\!3/2}\sum_{r=0}^\infty (-)^r C_r \bl(\frac{\la}{\pi^2a^2}\br)^{\!r} \sum_{k\geq 1}\frac{e^{-\pi^2k^2a^2/\la}}{k^{2r+2}}
\ee
as $|a|\to\infty$ in $|\arg\,a|<\pi/4$, where from (\ref{e37}) the coefficients $C_r$ are given by
\[C_r=\frac{2\g(r+\f{3}{2})}{\sqrt{\pi}}\,{}_1F_1(-r;\f{3}{2};\la).\]
The first few $C_r$ are therefore
\[C_0=1,\quad C_1=\f{3}{2}-\la,\quad C_2=\f{15}{4}-5\la+\la^2,\]
\[C_3=\f{105}{8}-\f{105}{4}\la+\f{21}{2}\la^2-\la^3.\]

We remark that the case $\mu=1$, $\ga=2$ can be obtained directly from Theorem 2, but also follows from the Poisson-Jacobi formula (\ref{e12a}), together with (\ref{e44}) and (\ref{e45}), since
\begin{eqnarray*}
S_{1,2}(a;\la)&=&\sum_{n=1}^\infty \frac{n^2 e^{-\la n^2/a^2}}{n^2+a^2}=\sum_{n=1}^\infty \bl(1-\frac{a^2}{n^2+a^2}\br) e^{-\la n^2/a^2}\\
&=&S_{0,0}(a;\la)-a^2 S_{1,0}(a;\la).
\end{eqnarray*}
\bigskip

\noindent {\bf Example 3.}\ \ Let $\mu=2$ and $\ga=0$ (so that $\delta=4$). Then we have
\bee\label{e46}
S_{2,0}(a;\la)=\sum_{n=1}^\infty \frac{e^{-\la n^2/a^2}}{(n^2+a^2)^2}=\frac{\sqrt{\pi}}{2a^3}\,U(\fs,-\fs,\la)-\frac{1}{2a^4}+J(a;\la).
\ee
From (\ref{e33}) and (\ref{e33c}), the coefficients $A_0=\la$, $A_1=-1$, $B_{00}=1$, $B_{01}=B_{11}=\fs$, so that
\[J_1(a;\la)=-\frac{\pi e^{\la-\pi a}}{2a^3 \sinh \pi a}\bl\{\la-\frac{1}{2}-\frac{\pi a e^{\pi a}}{2\sinh \pi a}\br\}\]
and
\[J_2(a;\la)\sim\bl(\frac{\la}{\pi a^2}\br)^{\!7/2} \sum_{r=0}^\infty (-)^r C_r \bl(\frac{\la}{\pi^2a^2}\br)^{\!r} \sum_{k\geq 1}\frac{e^{-\pi^2k^2a^2/\la}}{k^{2r+4}}.\]
From (\ref{e37}) the coefficients $C_r$ are given by
\[C_r=\frac{4(r+1)}{3\sqrt{\pi}}\,\g(r+\f{5}{2})\,{}_1F_1(-r;\f{5}{2};\la)\]
so that the first few coefficients are therefore
\[C_0=1,\quad C_1=5-2\la, \quad C_2=\f{105}{4}-21\la+3\la^2,\]
\[C_3=\f{315}{2}-189\la+54\la^2-4\la^3.\]
Then the exponentially small contribution is
\[J(a;\la)+\frac{\pi e^{\la-\pi a}}{2a^3 \sinh \pi a}\bl\{\la-\frac{1}{2}-\frac{\pi a e^{\pi a}}{2\sinh \pi a}\br\}\hspace{5cm}\]
\bee\label{e47}
\hspace{4cm}\sim \bl(\frac{\la}{\pi a^2}\br)^{\!7/2} \sum_{r=0}^\infty (-)^r C_r \bl(\frac{\la}{\pi^2a^2}\br)^{\!r} \sum_{k\geq 1}\frac{e^{-\pi^2k^2a^2/\la}}{k^{2r+4}}
\ee
as $|a|\to\infty$ in $|\arg\,a|<\pi/4$.

As mentioned in the previous example, the sums with $\mu=2$ and $p=2, 4$ can be obtained directly from Theorem 2, but also from the identities
\begin{eqnarray*}
S_{2,1}(a;\la)&=&S_{1,0}(a;\la)-a^2 S_{2,0}(a;\la)\\
S_{2,2}(a;\la)&=&S_{0,0}(a;\la)-2a^2 S_{2,1}(a;\la)-3a^4 S_{2,0}(a;\la).
\end{eqnarray*}
\vspace{0.6cm}

\begin{center}
{\bf 6.\ Numerical results and concluding remarks}
\end{center}
\setcounter{section}{6}
\setcounter{equation}{0}
\renewcommand{\theequation}{\arabic{section}.\arabic{equation}}
The expansion of the exponentially small contribution (when $\ga=2p$) given in Theorem 1 is exact for $\mu\geq 0$. It is possible to employ an asymptotic expansion for the integrals $I_{jk}$, but this would necessarliy introduce an error. However, we have evaluated these integrals to high numerical precision and have thereby verified the expansion (\ref{e322}) of $J(a;\la)$ for several parameter values to 50 decimal precision.

We present some numerical examples of the large-$a$ expansion of $S_{\mu,\ga}(a;\la)$ given in Theorem 2 to demonstrate the accuracy of our results. We subtract from the sum $S_{\mu,\ga}(a;\la)$ the finite terms appearing in (\ref{e38a}) by defining
\bee\label{e51}
{\hat S}_{\mu,\ga}(a;\la):=S_{\mu,\ga}(a;\la)-\{a^{1-\delta} {\cal H}(1)+H_{\mu,\ga}(a;\la)+J_1(a;\la)\}
\ee
and comparing it with the exponentially small asymptotic expansion $J_2(a;\la)$ in (\ref{e3j}). We stress that the contribution $J_1(a;\la)$ in (\ref{e34}) is an exact result when $\mu$ is an integer. In Table 1 we show the values of the absolute relative error in the high-precision computation of ${\hat S}_{\mu,\ga}(a;\la)$ from (\ref{e11}) using the asymptotic expansion for $J_2(a;\la)$ for different truncation index $r$. The values of $\mu$ and $\ga$ chosen correspond to the examples given in Section 5. The final entry in each column gives the value of ${\hat S}_{\mu,\ga}(a;\la)$. It is seen that the exponentially small contribution to $S_{\mu,\ga}(a;\la)$ when $\ga$ is an even integer agrees well with the expansion given in Theorem 2.
\begin{table}[t]
\caption{\footnotesize{The absolute relative error in the computation of ${\hat S}_{\mu,\ga}(a;\la)$ from (\ref{e51}) for different $\mu$, $\gamma$ and truncation index $r$ in the asymptotic expansion $J_2(a;\la)$ when $\lambda=2$ and $a=3$.}}
\begin{center}
\begin{tabular}{|l|c|c|c|}
\hline
&&&\\[-0.3cm]
\mcol{1}{|c|}{$r$} & \mcol{1}{c|}{$\mu=0,\ \ga=-2$} &\mcol{1}{c|}{$\mu=1,\ \ga=0$}& \mcol{1}{c|}{$\mu=2,\ \ga=0$} \\
[.1cm]\hline
&&&\\[-0.25cm]
0 & $3.307\times 10^{-02}$ & $1.007\times 10^{-02}$ & $2.464\times 10^{-02}$ \\
1 & $1.823\times 10^{-03}$ & $1.070\times 10^{-03}$ & $1.572\times 10^{-03}$ \\
2 & $1.408\times 10^{-04}$ & $5.900\times 10^{-05}$ & $3.764\times 10^{-04}$ \\
5 & $2.438\times 10^{-07}$ & $7.124\times 10^{-08}$ & $3.325\times 10^{-07}$ \\
10& $9.421\times 10^{-11}$ & $3.101\times 10^{-12}$ & $7.198\times 10^{-10}$ \\
15& $3.138\times 10^{-13}$ & $6.596\times 10^{-14}$ & $2.317\times 10^{-12}$ \\
20& $4.678\times 10^{-15}$ & $4.335\times 10^{-16}$ & $5.392\times 10^{-14}$ \\
[.1cm]\hline
&&&\\[-0.25cm]
${\hat S}_{\mu,\ga}$ & $-9.3737097\times 10^{-22}$ & $-9.7822227\times 10^{-22}$ & $+4.7287147\times 10^{-24}$\\
 [.2cm] \hline
\end{tabular}
\end{center}
\end{table}

It is worth mentioning that $J(a;\la)$ given in Theorems 1 and 2 appears to comprise two different types of exponentially small terms, namely $\exp (-2\pi ka)$ in Theorem 1 and both $\exp (-2\pi ka)$ and $\exp (-\pi^2 k^2a^2/\la)$, $k\geq 1$ in Theorem 2. However, a closer examination of the integrals $I_{jk}$ appearing in Theorem 1 reveals that they also contain the more subdominant terms $\exp (-\pi^2 k^2a^2/\la)$. To see this we consider
\[e^{-2\pi ka}I_{jk}=e^{-2\pi ka}\int_0^\infty \frac{t^{\mu-1}e^{-\psi(t)}}{(1+t)^\beta}dt,\qquad \beta=2p-j+\fs.\]
The phase function $\psi(t)$ in (\ref{e320}) has a saddle point at $t=t_s$, where $1+t_s=\sqrt{X_k/\la}=\pi ka/\la$ and $\psi(t_s)=0$, $\psi''(t_s)=2X_k/(1+t_s)^3=2\la^2/(\pi ka)$. For large complex $a$ in the sector $|\arg\,a|<\pi/4$, the integration path is chosen to emanate from the origin in the direction $\arg\,t=\phi=\pi-2\arg\,a$ to the singularity at $t=-1$ and thence along the path of steepest descent through $t_s$ to infinity in $\Re (t)>0$. The contribution from the saddle is controlled by
\[2e^{-2\pi ka} \sqrt{\frac{\pi}{2\psi''(t_s)}}\,\frac{t_s^{\mu-1}}{(1+t_s)^\beta}=\sqrt{\frac{\pi}{\la}} \bl(\frac{\pi ka}{\la}\br)^{\!\mu-\beta-1/2} e^{-2\pi ka}\]
while that from the  neighbourhood of the origin is approximately
\[e^{-X_k+i\mu\phi}\int_0^\infty e^{-|X_k|\tau}\tau^{\mu-1}d\tau=O(X_k^{-\mu} e^{-X_k}),\]
which produces the more subdominant exponential terms.

Finally, we note that the alternating version of (\ref{e11}) can be expressed in terms of $S_{\mu,\ga}(a;\la)$ since 
\[\sum_{n=1}^\infty \frac{(-)^{n-1}n^\ga}{(n^2+a^2)^\mu}\,e^{-\la n^2/a^2}=S_{\mu,\ga}(a;\la)-2^{1-\delta} S_{\mu,\ga}(\fs a;\la).\]
Application of Theorems 1 and 2 then enables the large-$a$ expansion of the alternating series to be determined.
\vspace{0.6cm}

\begin{center}
{\bf Appendix A: The pole structure of ${\cal H}(s)$}
\end{center}
\setcounter{section}{1}
\setcounter{equation}{0}
\renewcommand{\theequation}{\Alph{section}.\arabic{equation}}
The function ${\cal H}(s)$ defined in (\ref{e21b}), (\ref{e21c}) and (\ref{e21d}) has pokes at $s=-2k-\ga$, $k=0, 1, 2, \ldots\ $ and {\it apparent\/} poles at $s=\pm 2k+\delta$, $\delta=2\mu-\ga$. We shall show in this appendix that ${\cal H}(s)$ is regular at these last points. We have
\[{\cal H}(s)=\frac{\pi G(s)}{2\sin \pi(\mu-\f{\ga+s}{2})},\]
where
\[G(s):=\frac{\g(\f{\ga+s}{2})}{\g(\mu)} {\bf F}(\f{\ga+s}{2};1\!+\!\f{\ga+s}{2}\!-\!\mu;\la)-\la^{\mu-(\ga+s)/2} {\bf F}(\mu;1\!-\!\f{\ga+s}{2}\!+\!\mu;\la).\]
Here ${\bf F}$ denotes the normalised confluent hypergeometric function defined by  
\[{\bf F}(a;b,z)=\frac{1}{\g(b)} \,{}_1F_1(a;b;z),\]
which is defined for all values of the parameter $b$.                                               

Let $s_k=2k+\delta$ so that $\fs(\ga+s_k)=k+\mu$. Then
\begin{eqnarray*}
G(s_k)&=&(\mu)_k {\bf F}(\mu+k;k+1;\la)-\la^{-k} \sum_{r=k}^\infty \frac{(\mu)_r \la^r}{r! \g(1+r-k)}\\
&=&(\mu)_k {\bf F}(\mu+k;k+1;\la)-\sum_{r=0}^\infty \frac{(\mu)_{r+k} \la^r}{r! \g(1+k+r)}\\
&=&(\mu)_k {\bf F}(\mu+k;k+1;\la)-(\mu)_k \sum_{r=0}^\infty \frac{(\mu+k)_r \la^r}{r! \g(1+k+r)}\equiv 0.
\end{eqnarray*}
Hence ${\cal H}(s)$ is regular at $s_k=2k+\delta$.

A similar argument when $s_k=-2k+\delta$ shows that
\begin{eqnarray*}
G(s_k)&=&\frac{\g(\mu-k)}{\g(\mu)} \sum_{r=k}^\infty \frac{(\mu-k)_r \la^r}{r! \g(1-k+r)}-\la^k {\bf F}(\mu;1+k;\la)\\
&=&\frac{\g(\mu-k)}{\g(\mu)} \sum_{r=0}^\infty \frac{(\mu-k)_{r+k} \la^{r+k}}{r! \g(1+k+r)}-\la^k {\bf F}(\mu;1+k;\la)\\
&=&\la^k\sum_{r=0}^\infty \frac{(\mu)_r \la^r}{r! \g(1+k+r)}-\la^k {\bf F}(\mu;1+k;\la)\equiv 0,
\end{eqnarray*}
so that ${\cal H}(s)$ is also regular at the points $s_k=-2k+\delta$.

\vspace{0.6cm}

\begin{center}
{\bf Appendix B: The expansion in the case $\ga=-1$}
\end{center}
\setcounter{section}{2}
\setcounter{equation}{0}
\renewcommand{\theequation}{\Alph{section}.\arabic{equation}}
We consider the large-$a$ expansion of $S_{\mu,\ga}(a;\la)$ given in (\ref{e24}) in the special case $\ga=-1$ when the singularity of the integrand in (\ref{e22}) at $s=1$ is a double pole.
We set $s=1+\eps$, with $\eps\to0$. Then
\[{\cal H}_1(s)=\frac{\g(\fs\eps)\g(\mu-\fs\eps)}{2\g(\mu)}\,{}_1F_1(\fs\eps;1-\mu+\fs\eps;\la)\]
\[=\frac{1}{\eps}\bl\{1+\fs\eps(\psi(1)-\psi(\mu))+O(\eps^2)\br\}\,{}_1F_1(\fs\eps;1-\mu+\fs\eps;\la),\]
where
\[{}_1F_1(\fs\eps;1-\mu+\fs\eps;\la)=1+\frac{\eps\la}{2(1-\mu)}\bl\{1+\frac{1!\la}{(2-\mu)2!}+\frac{2!\la^2}{(2-\mu)_2 3!}+\cdots\br\}+O(\eps^2)\]
\[=1+\frac{\eps\la}{2(1-\mu)}\,{}_2F_2(1,1;2,2-\mu;\la)+O(\eps^2).\]
Using the fact that $\zeta(1+\eps)=\eps^{-1}(1+\ga_{_E}\eps+O(\eps^2))$ and $\psi(1)=-\ga_{_E}$, where $\ga_{_E}$ is the Euler-Mascheroni constant,
we obtain the residue  resulting from ${\cal H}_1(s)$ at $s=1$ given by
\[a\bl\{\log\,a+\frac{1}{2}\ga_{_E}
-\frac{1}{2}\psi(\mu)+\frac{\la}{2(1-\mu)}\,{}_2F_2(1,1;2,2-\mu;\la)\br\}\qquad(\mu\neq 1, 2, \ldots).\]
The residue resulting from ${\cal H}_2(s)$ is
\[a{\cal H}_2(1)=\frac{1}{2}\la^\mu \g(1-\mu)\,{}_1F_1(\mu;1+\mu;\la)\qquad (\mu\neq 1, 2, \ldots).\]

Hence, provided $\mu\neq 1, 2, \ldots\ $,
\[S_{\mu,-1}(a;\la)=\sum_{n=1}^\infty \frac{e^{-\la n^2/a^2}}{n(n^2+a^2)^\mu}\sim
a^{1-\delta}\bl\{\log\,a+\frac{1}{2}\ga_{_E}
-\frac{1}{2}\psi(\mu)+\frac{\la}{2(1-\mu)}\,{}_2F_2(1,1;2,2-\mu;\la)\br\}\]
\bee
+a^{1-\delta}\sum_{k=1}^\infty \frac{(-)^k (\mu)_k}{k!} \zeta(1-2k)\,{}_1F_1(-k;1-\mu-k;\la) a^{-2k}
\ee
as $a\to\infty$ in $|\arg\,a|<\pi/4$.

When $\mu$ is a positive integer a limiting process is required. To illustrate, we consider only the case $\mu=1$.
We find that ${\cal H}_1(1+\eps)=e^\la/\eps$ and
\[{\cal H}_2(1+\eps)=-\frac{\la}{\eps(1-\fs\eps)} \g(1+\fs\eps)\,{}_1F_1(1;2-\fs\eps;\la)\]
\[=-\frac{\la}{\eps}\bl\{1+\frac{1}{2}\eps(1-\ga_{_E}-\log\,\la)+O(\eps^2)\br\}\,{}_1F_1(1;2-\fs\eps;\la),\]
where
\begin{eqnarray*}
{}_1F_1(1;2-\fs\eps;\la)\!\!\!&=&\!\!\!1\!+\!\frac{\la}{2}\bl(1\!+\!\frac{\eps}{4}\br)\!+\!\frac{\la^2}{2\cdot 3}\bl(1\!+\!\frac{\eps}{4}\!+\!\frac{\eps}{6}\br)+\cdots +O(\eps^2)\\
&=&\!\!\!{}_1F_1(1;2;\la)+\frac{\eps}{2}\bl\{\frac{\la}{2} \frac{1}{2}\!+\!\frac{\la^2}{2\cdot 3}\bl(\frac{1}{2}\!+\!\frac{1}{3}\br)\!+\!\frac{\la^3}{2\cdot 3\cdot 4}\bl(\frac{1}{2}\!+\!\frac{1}{3}\!+\!\frac{1}{4}\br)\!+\!\cdots\br\}+O(\eps^2)\\
&=&\!\!\!\frac{e^\la-1}{\la}\!+\!\frac{\eps}{2}\sum_{n=1}^\infty \frac{\la^n \tau(n)}{(2)_n}+O(\eps^2),\qquad \tau(n):=\sum_{r=1}^n \frac{1}{r+1}.
\end{eqnarray*}

Then we obtain the expansion
\[S_{1,-1}(a;\la)=\sum_{n=1}^\infty \frac{e^{-\la n^2/a^2}}{n(n^2+a^2)}\sim
\frac{a^{1-\delta}}{2}\bl\{e^\la(\log\,\la+\ga_{_E}-1)+\log\,(a^2/\la)+\ga_{_E}+1-\la\sum_{n=1}^\infty \frac{\la^n \tau(n)}{(2)_n}\br\}\]
\bee
+a^{1-\delta} \sum_{k=1}^\infty \frac{(-)^k (\mu)_k}{k!} \zeta(1-2k) e_k(\la) a^{-2k}\ee
as $a\to\infty$ in $|\arg\,a|<\pi/4$, where
\[e_k(\la):={}_1F_1(-k;-k;\la)=\sum_{n=0}^k \frac{\la^n}{n!}.\]

\end{document}